\documentclass[11pt]{article}

\parskip=\medskipamount
\def\fpd#1#2{{\displaystyle\frac{\partial #1}{\partial #2}}}
\def\spd#1#2#3{{\displaystyle\frac{\partial^2 #1}
{\partial #2\partial #3}}}
\def\lie#1{{\cal L}_{#1}}
\def\vf#1{\frac{\partial}{\partial #1}}

\def\clift#1{#1^{\scriptscriptstyle{\mathrm{C}}}}

\def\vlift#1{#1^{\scriptscriptstyle{\mathrm{V}}}}

\def\E{\mathcal{E}}
\def\R{\mathcal{R}}

\def\onehalf{{\textstyle\frac12}}

\font\frak=eufm10 scaled\magstep1
\def\goth #1{\hbox{{\frak #1}}}
\def\g{\goth{g}}

\def\psiM{\psi^{\scriptscriptstyle M}}

\def\Ad{\mathop{\mathrm{ad}}\nolimits}
\def\ad{\Ad}

\begin{document}

\title{Relative equilibria of Lagrangian systems with symmetry}

\author{M.\ Crampin${}^{a}$ and T.\ Mestdag${}^{a,b}$\\[2mm]
{\small ${}^a$Department of Mathematical Physics and Astronomy, Ghent University}\\
{\small Krijgslaan 281, B-9000 Ghent, Belgium}\\[1mm]
{\small ${}^b$ Department of Mathematics, University of Michigan}\\
{\small 530 Church Street, Ann Arbor, MI 48109, USA}}

\date{}

\maketitle

{\small {\bf Abstract.} We discuss the characterization of relative
equilibria of Lagrangian systems with symmetry.
\\[2mm]
{\bf Mathematics Subject Classification (2000).}
34A05, 34A26, 37J15, 37J15, 70H03.
\\[2mm]
{\bf Keywords.} Lagrangian system, symmetry, relative equilibrium,
momentum map.}

\section{Introduction}

A relative equilibrium of a Lagrangian system which is invariant under
a Lie group $G$ is a motion of the system which coincides with an
orbit of a 1-parameter subgroup of $G$.  Relative equilibria are of
interest in particle dynamics \cite{Lawson,Lewis,Simo} of course, and
also in Riemannian and Finsler geometry~\cite{HG,Latifi,Szenthe1},
though there they are often studied under different names.

In this paper we consider Lagrangians which are arbitrary apart from
being subject to certain regularity conditions, with symmetry groups
which are also arbitrary except for being required to act freely and
properly on the space.  We prove a very general criterion for finding
relative equilibria:\ a relative equilibrium is a critical point of
the restriction of the energy to a level set of the momentum.  We
discuss the relation between this result and a different criterion for
relative equilibria given by Lewis \cite{Lewis}.  We also consider in
some detail two special cases, namely the case where the configuration
space is a Lie group \cite{HG,Latifi,Szenthe} and the case of a simple
mechanical system \cite{MRS}, and in the context of the latter we make
some remarks about the so-called Saari conjecture~\cite{Lawson}.  One
of our purposes in this paper, indeed, is to provide a single
framework for a variety of results about the conditions for the
existence of relative equilibria both in general and in particular
circumstances.

We shall use methods based on the consideration of frames adapted to
the group action, and velocity variables associated with such
frames, variables which are sometimes called quasi-velocities.  To
the best of our knowledge the study of relative equilibria by such
methods has not been carried out before, at least in recent times.
We have already used these methods in studying other aspects of
dynamical systems with symmetry \cite{Paper1,Paper2,MCPaper}, and
some derivations which are passed over rather quickly here are dealt
with at somewhat greater length in these references; nevertheless
the present paper is designed to be reasonably self-contained.

The basic relevant facts about group actions are discussed in Section
2. Section 3 is devoted to explaining our approach to Lagrangian
theory. The main result is proved in Section 4. In Section 5 the
alternative criterion for the existence of relative equilibria due to
Lewis is derived using our formalism. The applications
are discussed in Section~6.

\section{Preliminaries}

Suppose that $\psiM: G\times M \to M$ is a free and proper
left action of a connected Lie group $G$ on
a manifold $M$. (In using left actions we follow the convention of
Marsden and Ratiu \cite{MRbook,MRS}. Other authors, including for
example  Kobayashi and Nomizu \cite{KN}, use right actions; as a
consequence our formulae may differ in sign from those to be found
elsewhere in the literature.) The manifold $M$ is therefore a
principal fibre bundle with group $G$, over a base manifold $B$ say.
Let $\g$ be the Lie algebra of $G$. For any $\xi\in\g$, $\tilde\xi$
will denote the corresponding fundamental vector field on $M$, that
is, the infinitesimal generator of the 1-parameter group
$\psiM_{\exp(t\xi)}$ of transformations of $M$. Since $G$ is
connected, a tensor field on $M$ is $G$-invariant if and only if its
Lie derivatives by all fundamental vector fields vanish.  In
particular, a vector field $X$ on $M$ is invariant if
$[\tilde{\xi},X]=0$ for all $\xi\in\g$; indeed, it is sufficient
that $[\tilde{E}_a,X]=0$, $a=1,2,\ldots,\dim(\g)$, where $\{E_a\}$
is any basis of $\g$.

We will work with a (local) basis $\{X_i,\tilde{E}_a\}$ of vector
fields on $M$ adapted to the bundle structure, where the $\tilde{E}_a$
are fundamental vector fields corresponding to a basis of $\g$, and
the $X_i$, $i=1,2,\ldots,\dim(B)$, are $G$-invariant.  To obtain such
invariant vector fields we may introduce a principal connection on $M$
and a local basis of vector fields on $B$ (a coordinate basis for
example), and take for the $X_i$ the horizontal lifts to $M$ of these
vector fields, relative to the connection.  We call such a basis
$\{X_i,\tilde{E}_a\}$ a standard basis.  The pairwise brackets of the
elements of a standard basis are
\[
[X_i,X_j] = R^a_{ij} {\tilde E}_a,\quad
[X_i,{\tilde E}_a]=0,\quad\mbox{and}\quad
[{\tilde E}_a,{\tilde E}_b]=-C_{ab}^c{\tilde E}_c:
\]
the $R^a_{ij}$ are the components of the curvature of the connection,
regarded as a $\g$-valued tensor field, and the $C_{ab}^c$ are the
structure constants of $\g$ with respect to the basis $\{E_a\}$ (the
minus sign occurs because the fundamental vector fields behave as
right, not left, invariant vector fields on $G$).

Since we will be concerned with Lagrangian functions and their
corresponding Euler-Lagrange equations we must consider also certain
geometrical structures on the tangent bundle of $M$, which will be
denoted by $\tau:TM\to M$.  One important idea is that of lifting
vector fields from $M$ to $TM$.  There are in fact two canonical
ways of carrying this out (see for example \cite{CP,YI} for more
details on the following material).  Let $Z$ be a vector field on
$M$.  The complete or tangent lift of $Z$ to $TM$, $\clift{Z}$, is
the vector field whose flow consists of the tangent maps of the flow
of $Z$.  The vertical lift of $Z$, $\vlift{Z}$, is tangent to the
fibres of $\tau$ and on the fibre over $m\in M$ coincides with the
constant vector field $Z_m$.  Then $T\tau(\clift{Z})=Z$ while
$T\tau(\vlift{Z})=0$. Complete and vertical lifts
satisfy the following bracket relations:
\[
[\clift{Y},\clift{Z}]=\clift{[Y,Z]},\quad
[\clift{Y},\vlift{Z}]=\vlift{[Y,Z]},\quad\mbox{and}\quad
[\vlift{Y},\vlift{Z}]=0.
\]
From a standard basis $\{X_i,{\tilde E}_a\}$ on $M$ we may construct
a standard basis $\{\clift{X_i},\clift{\tilde{E}_a},
\vlift{X_i},\vlift{\tilde{E}_a}\}$ on $TM$ by taking complete and
vertical lifts. We will need to use the following bracket relations
satisfied by these vector fields:
\[
[\clift{{\tilde E}_a},\clift{X_i}]=
[\clift{{\tilde E}_a},\vlift{X_i}]=0,\quad
[\clift{{\tilde E}_a},\clift{{\tilde E}_b}]
=-C^c_{ab} \clift{{\tilde E}_{c}},\quad
[\clift{\tilde{E}_a},\vlift{\tilde{E}_b}]
=-C_{ab}^c\vlift{\tilde{E}_c}.
\]

We can use any basis of vector fields $\{Z_\alpha\}$ on a manifold $M$
to introduce fibre coordinates on $TM$, simply by taking the
coordinates of a point $u$ in the fibre over $m$ to be the components
of $u\in T_mM$ with respect to the basis $\{Z_\alpha|_m\}$ of $T_mM$;
such fibre coordinates are sometimes called quasi-velocities, and we
will follow this practice.  We can specify quasi-velocities more
succinctly as follows.  Let $\{\theta^\alpha\}$ be the basis of
1-forms on $M$ dual to the basis $\{Z_\alpha\}$ of vector fields, and
for any 1-form $\theta$ on $M$ let $\hat{\theta}$ denote the function
on $TM$ defined by $\hat{\theta}(m,u)=\langle u,\theta_m\rangle$.
Then the functions $\hat{\theta}^\alpha$ are the quasi-velocities
corresponing to the $Z_\alpha$.  The calculation of the derivatives of
quasi-velocities along complete and vertical lifts of basis vector
fields is carried out with the use of the following formulae:
\[
\clift{Z}(\hat{\theta})=\widehat{\lie{Z}\theta}, \quad
\vlift{Z}(\hat{\theta})=\tau^*\theta(Z).
\]
In particular,
$\vlift{Z_\alpha}(\hat{\theta}^\beta)=\delta^\beta_\alpha$.

Consider now a standard basis $\{X_i,\tilde{E}_a\}$.  We write
$(v^i,v^a)$ for the corresponding quasi-velocities.  Using the
formulae above, we obtain
\[
\begin{array}{lllllll}
\clift{X_i}(v^j)=0,&\quad&\vlift{X_i}(v^j)=\delta^j_i,&\quad&
\clift{X_i}(v^a)=-R^a_{ij}v^j,&\quad&\vlift{X_i}(v^a)=0,\\
\clift{\tilde{E}_a}(v^i)=0,&&\vlift{\tilde{E}_a}(v^i)=0,&&
\clift{\tilde{E}_a}(v^b)=C_{ac}^bv^c,&&\vlift{\tilde{E}_a}(v^b)=\delta^b_a.
\end{array}
\]
It will sometimes be convenient to use a slightly unconventional
notation for points in $TM$:\ we will denote such points in the form
$(m,v^i,v^a)$, where $(v^i,v^a)$ are the quasi-velocities of a point
in $T_mM$ with respect to a specific standard basis.

\section{The Euler-Lagrange equations}

We next explain our approach to Lagrangian theory, beginning with the
general situation where no symmetries are assumed.

A Lagrangian $L$ is a function on a tangent bundle $TM$ (we deal only
with the autonomous case).  Take local
coordinates $(x^\alpha)$ on $M$ and the corresponding local
coordinates $(x^\alpha,u^\alpha)$ on $TM$.  The Euler-Lagrange
equations of $L$,
\[
\frac{d}{dt}\left(\fpd{L}{u^\alpha}\right)-\fpd{L}{x^\alpha}=0,
\]
are second-order ordinary differential equations for
the extremals. However, the second derivatives $\ddot{x}^\alpha$
are not necessarily determined by these equations.
We say that $L$ is regular if
\[
\spd{L}{u^\alpha}{u^\beta},
\]
its Hessian with respect to the fibre coordinates, is everywhere
non-singular when considered as a symmetric matrix.  When the
Lagrangian is regular the Euler-Lagrange equations may be solved
explicitly for the $\ddot{x}^\alpha$ to give a system of differential
equations of the form $\ddot{x}^\alpha=\Gamma^\alpha(x,\dot{x})$; in
turn, these equations can be thought of as defining a vector field
$\Gamma$ on $TM$, namely
\[
\Gamma=u^\alpha\vf{x^\alpha}+\Gamma^\alpha\vf{u^\alpha}.
\]
This vector field, which is an example of a second-order differential
equation field, is called the Euler-Lagrange field of $L$.  The
Euler-Lagrange equations may be written
\[
\Gamma\left(\fpd{L}{u^\alpha}\right)-\fpd{L}{x^\alpha}=0;
\]
they determine $\Gamma$, assuming it to be a second-order differential
equation field, when $L$ is regular.

In this paper we will assume that $L$ is regular and we will work with
the Euler-Lagrange equations in terms of the second-order differential
equation field $\Gamma$.  However, we need to be able to express those
equations, and the property of being a second-order differential
equation field, in terms of a basis of vector fields on $M$ which is not
necessarily of coordinate type, say $\{X_\alpha\}$.
A vector field is a second-order differential equation field if it
takes the form
\[
\Gamma=\hat{u}^\alpha\clift{X_\alpha}+\hat{\Gamma}^\alpha\vlift{X_\alpha}
\]
where the $\hat{u}^\alpha$ are the quasi-velocities corresponding to
the basis $\{X_\alpha\}$.  Furthermore, the equations
\[
\Gamma(\vlift{X_\alpha}(L))-\clift{X_\alpha}(L)=0
\]
are equivalent to the Euler-Lagrange equations.

We will also need a coordinate-independent expression for the Hessian.
In fact the Hessian $g$ of $L$, evaluated at $u\in TM$, is the
symmetric bilinear form $g_u$ on $T_mM$, $m=\tau(u)$, given by
$g_u(v,w)=\vlift{v_u}(\vlift{w}(L))$, where the vertical lifts are
considered as vector fields on $T_mM$.  We can equally well regard
$g_u$ as a bilinear form on the vertical subspace of $T_uTM$, by
identifying $v$ and $w$ with their vertical lifts.  Since we assume
that $L$ is regular we know that $g$ is non-singular.

Suppose now that $L$ has a symmetry group $G$, acting to the left on
$M$ in such a way that $M$ is a principal bundle with $G$ as its
group, as we described above.  By saying that $G$ is a symmetry
group of the Lagrangian we mean that $L$ is invariant under the
induced action of $G$ on $TM$, so that $\clift{\tilde{\xi}}(L)=0$
for all $\xi\in\g$. A regular invariant Lagrangian determines an
Euler-Lagrange field which is also invariant.

We choose a standard basis of vector fields $\{X_i,\tilde{E}_a\}$ on
$M$, as described above.  The invariance of the Lagrangian can
be characterized by the property $\clift{\tilde E}_a(L)=0$.
The Euler-Lagrange equations for $L$ are
\begin{eqnarray*}
\Gamma(\vlift{X_i}(L))-\clift{X_i}(L)&=&0\\
\Gamma(\vlift{\tilde{E}_a}(L))-\clift{\tilde{E}_a}(L)&=&0.
\end{eqnarray*}
It follows immediately from invariance that
$\Gamma(\vlift{\tilde{E}_a}(L))=0$, which is to say that the functions
$\vlift{\tilde{E}_a}(L)$, which we denote by $p_a$, are first
integrals of $\Gamma$.  In fact the $p_a$ can be regarded as
components of an element of $\g^*$, the dual of the Lie algebra $\g$,
and the corresponding vector is called the momentum.  The map
$TM\to\g^*$ by $v\mapsto(p_a(v))$ is equivariant between the
given action of $G$ on $TM$ and the coadjoint action of $G$ on $\g^*$.
We have
\[
\clift{\tilde{E}_a}(p_b)=\clift{\tilde{E}_a}\vlift{\tilde{E}_b}(L)
=[\clift{\tilde{E}_a},\vlift{\tilde{E}_b}](L)
=-C_{ab}^c\vlift{\tilde{E}_c}(L)=-C_{ab}^cp_c,
\]
which expresses the differential version of this result in our formalism.

The Euler-Lagrange field $\Gamma$ is tangent to any level set of
momentum, that is, any subset of $TM$ of the form
$p_a=\mu_a=\mbox{constant}$, $a=1,2,\ldots,\dim G$ --- provided of
course that it is a submanifold.  To describe when this is so we
have recourse to the Hessian again.  The components of the Hessian
$g$ with respect to our standard basis will be expressed as follows:
\[
g(\tilde{E}_a,\tilde{E}_b)=g_{ab},\quad g(X_i,X_j)=g_{ij},\quad
g(X_i,\tilde{E}_a)=g_{ia}=g_{ai}=g(\tilde{E}_a,X_i)
\]
(in general these will be functions on $TM$, not $M$).  Then if
$(g_{ab})$ is non-singular the equations $p_a=\mu_a$ in principle
determine the $v_a$ in terms of the other variables, so the level set
of momentum will be a submanifold; we accordingly make the further
assumption about $L$ that $(g_{ab})$ is non-singular everywhere.

We will be working on a level set of momentum, say $p_a=\mu_a$, which
we denote by $N_\mu$. We will next define
vector fields related to $\clift{X_i}$, $\vlift{X_i}$ and
$\clift{\tilde{E}_a}$ which are tangent to $N_\mu$.
Since by assumption $(g_{ab})$ is non-singular, there are uniquely defined
coefficients $A^b_i$, $B^b_i$ and $C^b_a$ such that
\begin{eqnarray*}
(\clift{X_i}+A^b_i\vlift{\tilde{E}_b})(p_a)&=&
\clift{X_i}(p_a)+A^b_ig_{ab}=0\\
(\vlift{X_i}+B^b_i\vlift{\tilde{E}_b})(p_a)&=&
\vlift{X_i}(p_a)+B^b_ig_{ab}=0\\
(\clift{\tilde{E}_a}+C^b_a\vlift{\tilde{E}_b})(p_c)&=&
\clift{\tilde{E}_a}(p_c)+C^b_ag_{bc}=0.
\end{eqnarray*}
Define vector fields $\clift{\bar{X}_i}$, $\vlift{\bar{X}_i}$ and
$\clift{\bar{E}_a}$ by
\begin{eqnarray*}
\clift{\bar{X}_i}&=&\clift{X_i}+A^a_i\vlift{\tilde{E}_a}\\
\vlift{\bar{X}_i}&=&\vlift{X_i}+B^a_i\vlift{\tilde{E}_a}\\
\clift{\bar{E}_a}&=&\clift{\tilde{E}_a}+C^b_a\vlift{\tilde{E}_b};
\end{eqnarray*}
they are tangent to each level set $N_\mu$.  (The notation is not
meant to imply that the barred vector fields are actually complete or
vertical lifts.)  We need expressions for the actions of
$\clift{\bar{X}_i}$, $\vlift{\bar{X}_i}$ and $\clift{\bar{E}_a}$ on
$v^i$ and $v^a$, and for their pairwise brackets. For the former we
have
\[
\begin{array}{lll}
\vlift{\bar{X}_i}(v^j)=\delta^j_i,\quad&
\clift{\bar{X}_i}(v^j)=0,\quad&
\clift{\bar{E}_a}(v^i)=0,\\
\vlift{\bar{X}_i}(v^a)=B^a_i,\quad&
\clift{\bar{X}_i}(v^a)=-R^a_{ij}v^j+A^a_i,\quad&
\clift{\bar{E}_a}(v^b)=C_{ac}^bv^c+C^b_a.
\end{array}
\]
To find the brackets of barred vector fields we argue as follows.  The
vector fields $\vlift{\tilde{E}_a}$ are transverse to the level sets,
and the barred vector fields span them.  Thus on any level set the
bracket of any two of the barred vector fields is a linear combination
of vector fields of the same form.  Consider for example
$[\clift{\bar{E}_a},\vlift{\bar{X}_i}]$.  It is easy to see
from the expressions for $\clift{\bar{E}_a}$ and $\vlift{\bar{X}_i}$
that this bracket is at worst a linear combination of the
$\vlift{\tilde{E}_a}$; it follows immediately that
$[\clift{\bar{E}_a},\vlift{\bar{X}_i}]=0$.  By similar arguments we
can show that the brackets of the barred vector fields just
reproduce those of their unbarred counterparts, except that
$[\clift{\bar{X}_i},\vlift{\bar{X}_j}]=0$ (though we won't actually
use this fact).

We will now rewrite the Euler-Lagrange equations
$\Gamma(\vlift{X_i}(L))-\clift{X_i}(L)=0$, taking into account the
fact that $\Gamma$ is tangent to the level sets of momentum.  For this
purpose we introduce the function $\R$ on $TM$ given by
\[
\R=L-v^ap_a.
\]
Since $\R$ generalizes in an obvious way the classical Routhian
corresponding to ignorable coordinates \cite{MRbook,Routh} we call it
the {\em Routhian}.  We have discussed the generalization of Routh's
procedure to arbitrary regular Lagrangians with non-Abelian symmetry
groups elsewhere \cite{Paper2}; we must repeat the derivation of the
expression of the remaining Euler-Lagrange equations in terms of $\R$.

To obtain the desired equations we first express $\clift{X_i}(L)$ and
$\vlift{X_i}(L)$ in terms of the barred vector fields and the
Routhian, as follows:
\begin{eqnarray*}
\clift{X_i}(L)&=&\clift{\bar{X}_i}(L)-A^a_i\vlift{\tilde{E}_a}(L)\\
&=&\clift{\bar{X}_i}(L-v^ap_a)+(-R^a_{ij}v^j+A^a_i)p_a
+v^a\clift{\bar{X}_i}(p_a)-A^a_ip_a\\
&=&\clift{\bar{X}_i}(\R)-p_aR^a_{ij}v^j;\\
\vlift{X_i}(L)&=&\vlift{\bar{X}_i}(L)-B^a_i\vlift{\tilde{E}_a}(L)\\
&=&\vlift{\bar{X}_i}(L-v^ap_a)+B^a_ip_a
+v^a\vlift{\bar{X}_i}(p_a)-B^a_ip_a\\
&=&\vlift{\bar{X}_i}(\R).
\end{eqnarray*}
Thus if we denote by $\R^\mu$ the restriction of the Routhian to the
submanifold $N_\mu$ (where it becomes $L-v^a\mu_a$), taking account of
the fact that $\Gamma$ is tangent to $N_\mu$ we have
\[
\Gamma(\vlift{\bar{X}_i}(\R^\mu))-\clift{\bar{X}_i}(\R^\mu)=-\mu_aR^a_{ij}v^j.
\]
These are the reduced Euler-Lagrange equations, or the generalized Routh
equations as they are called in \cite{Paper2}.

Since $\Gamma$ satisfies $\Gamma(p_a)=0$ it may be expressed in the
form
\[
\Gamma=v^i\clift{\bar{X}_i}+\Gamma^i\vlift{\bar{X}_i}+v^a\clift{\bar{E}_a}.
\]
If the matrix-valued function $\vlift{\bar X}_i(\vlift{\bar X}_j(\R))$
is non-singular, the generalized Routh equations will determine the
coefficients $\Gamma^i$.  We show now that this is the case, as always
under the assumptions that $L$ is regular and that $(g_{ab})$ is
non-singular.

Recall that $\vlift{\bar{X}_i}=\vlift{X_i}+B^a_i\vlift{\tilde{E}_a}$
is determined by the condition that $\vlift{\bar{X}_i}(p_a)=0$; it
follows that $B^a_i=-g^{ab}g_{ib}$, where $(g^{ab})$ is the matrix
inverse to $(g_{ab})$.  Now $\vlift{\bar{X}_i}(\R)=\vlift{X_i}(L)$,
so
\[
\vlift{\bar{X}_i}(\vlift{\bar{X}_j}(\R))=
(\vlift{X_i}-g^{ab}g_{ib}\vlift{\tilde{E}_a})(\vlift{X_j}(L))
=g_{ij}-g^{ab}g_{ia}g_{jb}.
\]
It is a straightforward exercise in linear algebra to show that under the
stated conditions the matrix with these components is non-singular.

\section{Relative equilibria}

Consider an autonomous second-order differential equation field
$\Gamma$ on the tangent bundle $TM$ of a manifold $M$.  Let
$t\mapsto\gamma(t)$ be a base integral curve of $\Gamma$, that is, a
curve on $M$ whose natural lift $t\mapsto(\gamma(t),\dot{\gamma}(t))$
to $TM$ is an integral curve of $\Gamma$.  The curve $\gamma$ is
uniquely determined by its initial conditions
$(\gamma(0),\dot{\gamma}(0))$ and the fact that it is a base integral
curve.

Now suppose that a Lie group $G$ acts to the left on $M$ in such a
way that $M$ is a principal $G$-bundle, $\pi:M\to B$; and suppose
that $\Gamma$ is invariant under the induced action of $G$ on $TM$.
Then $G$ maps base integral curves of $\Gamma$ to base integral
curves; and for $g\in G$, $t\mapsto \psi^M_g(\gamma(t))$ is the base integral
curve with initial conditions
$\psi^{TM}_g(\gamma(0)),\dot{\gamma}(0))$.

A base integral curve $\gamma$ is a {\em relative equilibrium\/} of
$\Gamma$ if it coincides with an integral curve of a fundamental
vector field of the action of $G$ on $M$, that is, if
$\gamma(t)=\psi^M_{\exp(t\xi)}(m)$ for some $m\in M$, $\xi\in\g$; of course
$m=\gamma(0)$, and $\dot{\gamma}(0)=\tilde{\xi}_m$.  A relative
equilibrium is a curve in a fibre of $\pi:M\to B$, so that
$\pi(\gamma(t))$ is a fixed point of $B$; but not all curves that
project onto fixed points of $B$ are relative equilibria.  Evidently
if $\gamma$ is a relative equilibrium, so is $\psi^M_g\circ\gamma$ for any
$g\in G$.

The base integral curve $\gamma$ is a relative equilibrium if and
only if its natural lift coincides with an integral curve of a
fundamental vector field of the induced action of $G$ on $TM$.  That
is to say, if an integral curve of the vector field $\Gamma$
coincides with an integral curve of $\clift{\tilde{\xi}}$ for some
$\xi\in\g$ then the corresponding base integral curve is a relative
equilibrium, and conversely.  But
we are now dealing directly with an invariant vector field, namely
$\Gamma$; by invariance, the integral curve of $\Gamma$ through
$v\in TM$ will coincide with that of $\clift{\tilde{\xi}}$ if and
only if $\Gamma_v=\clift{\tilde{\xi}}_v$.  Thus finding relative
equilibria is a matter of locating points $v\in TM$ with the
property that $\Gamma_v=\clift{\tilde{\xi}}_v$ for some $\xi\in\g$;
we call such points relative equilibrium points.  We will shortly
address the problem of finding relative equilibrium points for the
Euler-Lagrange field of an invariant Lagrangian.

Recall that in the absence of symmetry, the equilibrium points of a
regular Lagrangian --- the zeros of its Euler-Lagrange field --- are
just the critical points of the energy.  It may be worth seeing why,
for comparison with what follows.  Let $(x^\alpha,u^\alpha)$ denote coordinates
on $TM$.  If $\E$ is the energy of a Lagrangian $L$, so that
\[
\E=u^\beta\fpd{L}{u^\beta}-L,
\]
then
\begin{eqnarray*}
\fpd{\E}{x^\alpha}&=&u^\beta\spd{L}{x^\alpha}{u^\beta}-\fpd{L}{x^\alpha}
=-\Gamma^\beta\spd{L}{u^\alpha}{u^\beta}+u^\beta\left(\spd{L}{x^\alpha}{u^\beta}
-\spd{L}{x^\beta}{u^\alpha}\right)\\
\fpd{\E}{u^\alpha}&=&u^\beta\spd{L}{u^\alpha}{u^\beta},
\end{eqnarray*}
and the critical points of $\E$ are precisely the points where
$u^\alpha=0$ and $\Gamma^\alpha=0$.

We will use these remarks as a guide to the formulation of a similar
result about relative equilibrium points in the Lagrangian formalism.
The energy $\E$ of the Lagrangian $L$ is given by $\E=\Delta(L)-L$,
where $\Delta$ is the Liouville field.  We note first that since
$[\Delta,\clift{Z}]=0$ for any vector field $Z$ on $M$, when $L$
is invariant $\E$ is also invariant.

We want an expression for the energy $\E$ of a Lagrangian in
terms of a standard basis, for which we need to know how to write
$\Delta$ with respect to such a basis:\ the obvious guess, namely
\[
\Delta=v^i\vlift{X_i}+v^a\vlift{\tilde{E}_a},
\]
is in fact correct. Thus the energy of $L$ is
\[
\E=\Delta(L)-L=v^i\vlift{X_i}(L)-(L-v^ap_a)=v^i\vlift{X_i}(L)-\R.
\]
We showed above that
$\vlift{X_i}(L)=\vlift{\bar{X}_i}(\R)$,
so we can write this as
\[
\E=v^i\vlift{\bar{X}_i}(\R)-\R.
\]

Next we derive expressions for the derivatives of $\E$ along the
barred vector fields. In the first place,
\[
\vlift{\bar{X}_i}(\E)=v^j\vlift{\bar{X}_i}(\vlift{\bar{X}_j}(\R)).
\]
Secondly, we have
$\clift{\bar{E}_a}(L)=
(\clift{\tilde{E}_a}+C_a^b\vlift{\tilde{E}_b})(L)=C_a^bp_b$,
so that
$\clift{\bar{E}_a}(\R)=C_a^bp_b-p_b(C_{ac}^bv^c+C_a^b)=
-C_{ac}^bp_bv^c$,
whence
\[
\clift{\bar{E}_a}(\E)=
\clift{\bar{E}_a}(v^i\vlift{\bar{X}_i}(\R)-\R)=
v^i\clift{\bar{E}_a}(\vlift{\bar{X}_i}(\R))-\clift{\bar{E}_a}(\R)=
C_{ac}^bp_bv^c+S_{ai}v^i,
\]
where $S_{ai}$ stands for an expression whose details will not concern
us. Finally,
\[
\clift{\bar{X}_i}(\E)=
\clift{\bar{X}_i}(v^j\vlift{\bar{X}_j}(\R)-\R)=
v^j\clift{\bar{X}_i}(\vlift{\bar{X}_j}(\R))-\clift{\bar{X}_i}(\R).
\]
From the generalized Routh equations
$\Gamma(\vlift{\bar{X}_i}(\R))-\clift{\bar{X}_i}(\R)=-\mu_aR^a_{ij}v^j$,
with $\Gamma$  expressed in the form
$\Gamma=v^i\clift{\bar{X}_i}+\Gamma^i\vlift{\bar{X}_i}+v^a\clift{\bar{E}_a}$,
we obtain
\[
\clift{\bar{X}_i}(\R)=\Gamma^j\vlift{\bar{X}_i}(\vlift{\bar{X}_j}(\R))
+v^a\clift{\bar{E}_a}\vlift{\bar{X}_i}(\R)+T_{ij}v^j
\]
where the exact form of $T_{ij}$ will again be of no concern. Now
\[
\clift{\bar{E}_a}(\vlift{\bar{X}_i}(\R))=
\vlift{\bar{X}_i}(\clift{\bar{E}_a}(\R))=
-\vlift{\bar{X}_i}(C_{ac}^bp_bv^c)=
-C_{ac}^bp_b\vlift{\bar{X}_i}(v^c).
\]
Thus
\[
\clift{\bar{X}_i}(\E)=-\Gamma^j\vlift{\bar{X}_i}(\vlift{\bar{X}_j}(\R))
+U_{ij}v^j+V^c_iv^aC_{ac}^bp_b,
\]
where $U_{ij}$ and $V^c_i$ will likewise be of no particular immediate
interest (though in fact $V^c_i=-B^c_i$).

We next consider the conditions for a relative equilibrium.  The
integral curve of $\Gamma$ through a point $(m,v^i,v^a)$ of $TM$ will
coincide with the integral curve of some $\clift{\tilde{\xi}}$,
$\xi\in\g$, if and only if
$\Gamma(m,v^i,v^a)=\clift{\tilde{\xi}}(m,v^i,v^a)$, that is, if and
only if at that point
\[
v^i\clift{\bar{X}_i}+\Gamma^i\vlift{\bar{X}_i}+v^a\clift{\bar{E}_a}
=\xi^a\clift{\tilde{E}_a}.
\]
Thus the integral curve of $\Gamma$ through a point $(m,v^i,v^a)$ of $TM$
will coincide with the integral curve of some $\clift{\tilde{\xi}}$
if and only if
\[
v^i=0,\quad v^a=\xi^a,\quad \Gamma^i(m,0,\xi^a)=0;
\]
and moreover we must have
$\xi^a\clift{\bar{E}_a}=\xi^a\clift{\tilde{E}_a}$, which just says
that the integral curve of $\clift{\tilde{\xi}}$ must lie in the level
set containing the point $(m,v^i,v^a)$, as does the integral curve of
$\Gamma$.  Let us assume that we are on the level set $p_a=\mu_a$;
then this last condition becomes
\[
\xi^aC^c_{ab}\mu_c=0.
\]
We can now prove that the relative equilibrium points lying in any
level set $N_\mu$ are just the critical points of
$\E^\mu$, the restriction of $\E$ to $N_\mu$, assuming as before that
$L$ is regular and $(g_{ab})$ is non-singular.

Suppose first there is a relative equilibrium point in $N_\mu$:\ it is a
point $(m,v^i,v^a)$ such that $v^i=0$, $v^aC^c_{ab}\mu_c=0$ and
$\Gamma^i(m,0,v^a)=0$.  From the formulae for the derivatives of $\E$
obtained above, we have $\vlift{\bar{X}_i}(\E^\mu)=
\clift{\bar{E}_a}(\E^\mu)=\clift{\bar{X}_i}(\E^\mu)=0$ at $(m,0,v^a)$,
and since these vector fields span the tangent distribution to the
level set, the point is a critical point of $\E^\mu$.

Conversely, suppose that a point $(m,v^i,v^a)$, lying in $N_\mu$, is a
critical point of $\E^\mu$, so that
$\vlift{\bar{X}_i}(\E^\mu)=\clift{\bar{E}_a}(\E^\mu)=\clift{\bar{X}_i}(\E^\mu)
=0$ there.  Since by assumption the symmetric-matrix-valued function
$\vlift{\bar{X}_i}(\vlift{\bar{X}_j}(\R))$ is non-singular, we find from
the condition $\vlift{\bar{X}_i}(\E^\mu)=0$ that $v^i=0$; from the
condition $\clift{\bar{E}_a}(\E^\mu)=0$ we obtain
$v^aC^c_{ab}\mu_c=0$; and from the condition
$\clift{\bar{X}_i}(\E^\mu)=0$ we deduce that $\Gamma^i(m,0,v^a)=0$.
The integral curve of $\Gamma$ through the point therefore coincides
with that of $\clift{\tilde{\xi}}$ where $\xi^a=v^a$.

As we have mentioned, the condition $\xi^aC^c_{ab}\mu_c=0$ states that
the fundamental vector field $\clift{\tilde{\xi}}$ is tangent to the
level set $N_\mu$. There is another way of interpreting this
condition.  We pointed out earlier that the map $v\mapsto(p_a(v))$
is equivariant between the given action of $G$ on $TM$ and the
coadjoint action of $G$ on $\g^*$.  For any $\mu\in\g^*$ we denote by
$G_\mu$ the isotropy group of $\mu$ under the coadjoint action, and
$\g_\mu$ its Lie algebra. By equivariance, $\xi\in\g_\mu$ if and only
if $\clift{\tilde{\xi}}$ is tangent to $N_\mu$. Thus
$\xi^aC^c_{ab}\mu_c=0$ is also the necessary and sufficient condition
that $\xi\in\g_\mu$.

\section{Lewis's criterion for a relative equilibrium}

We next discuss the somewhat different criterion for the existence of a
relative equilibrium given by Lewis in \cite{Lewis}.
Lewis defines the {\em locked Lagrangian\/} for any
$\xi\in\g$, $L_\xi$, by
\[
L_\xi(m)=L(m,\tilde{\xi}_m);
\]
thus $L_\xi$ is a function on $M$. She shows that a point
$(m,0,\xi^a)$ of $TM$ is a relative equilibrium point, for a regular
Lagrangian, if and only if $m$ is a critical point of $L_\xi$. We
now establish a similar result by our methods.

The first task is to relate the derivatives of $L_\xi$ to those of
$L$.  For this purpose it is helpful to observe that the
specification of $L_\xi$ can be regarded as a particular case of a
general construction.  Let $F$ be any function on the tangent bundle
$TM$ of some manifold $M$, and $X$ any vector field on $M$.  Then
$X$ is, or defines, a section of $TM\to M$, which we will denote by
$\sigma_X$ for clarity; and we can use such a section to obtain from
$F$ a function $F_X$ on $M$ by pull-back:\ $F_X=\sigma_X^*F$.  The
locked Lagrangian is an example of this construction, with $F=L$,
$X=\tilde{\xi}$.

We require a formula for $Y(F_X)$, the derivative of $F_X$ along any
other vector field $Y$ on $M$.  Now there is a unique vector field
$T\sigma_{X}(Y)$ on the image of the section $\sigma_X$ which is tangent
to it and which projects onto $Y$.  In fact for any $v\in T_mM$, say
$v^\alpha\partial/\partial x^\alpha$, the vector
\[
v^\alpha\vf{x^\alpha}+v^\beta\fpd{X^\alpha}{x^\beta}\vf{u^\alpha}
\in T_{\sigma_X(m)}TM
\]
is the unique vector which projects onto $v$ and is tangent to the
section. Thus
\[
T\sigma_{X}(Y)=Y^\alpha\vf{x^\alpha}+Y^\beta\fpd{X^\alpha}{x^\beta}\vf{u^\alpha}.
\]
Notice that we can express the right-hand side as
\[
Y^\alpha\vf{x^\alpha}+Y^\beta\fpd{X^\alpha}{x^\beta}\vf{u^\alpha}=
Y^\alpha\vf{x^\alpha}+X^\beta\fpd{Y^\alpha}{x^\beta}\vf{u^\alpha}
-\left(X^\beta\fpd{Y^\alpha}{x^\beta}
-Y^\beta\fpd{X^\alpha}{x^\beta}\right)\vf{u^\alpha},
\]
and this is just the restriction to the image of $\sigma_X$ of the
vector field $\clift{Y}-\vlift{[X,Y]}$, a vector field which is
defined globally on $TM$. Thus
\[
Y(F_X)=Y(\sigma_X^*F)=\sigma_X^*(T\sigma_{X}Y(F))=
\sigma_X^*\left((\clift{Y}-\vlift{[X,Y]})(F)\right).
\]

We now use this result to obtain expressions for the derivatives of
$L_\xi$ along the local basis vector fields $\tilde{E}_a$, $X_i$ on
$M$. We have $\clift{\tilde{E}_a}(L)=0$, while
$[\tilde{\xi},\tilde{E}_a]=C_{ab}^c\xi^b\tilde{E}_c$, whence
\[
\tilde{E}_a(L_\xi)=-C_{ab}^c\xi^b\sigma_{\tilde{\xi}}^*(p_c).
\]
On the other hand
\[
X_i(L_\xi)=\sigma_{\tilde{\xi}}^*(\clift{X_i}(L)),
\]
because $X_i$ is invariant under the $G$-action. But from the
Euler-Lagrange equations
$\clift{X_i}(L)=\Gamma(\vlift{X_i}(L))=\Gamma(\vlift{\bar{X}_i}(\R))$. So
finally, at any $m\in M$,
\begin{eqnarray*}
\tilde{E}_a|_m(L_\xi)&=&-C_{ab}^c\xi^b\mu_c\\
X_i|_m(L_\xi)&=&\Gamma(\vlift{\bar{X}_i}(\R))(m,0,\xi^a),
\end{eqnarray*}
where we have set $p_a(m,0,\xi^b)=\mu_a$.

Now suppose that $(m,0,\xi^a)$ is a relative equilibrium point on the
level set $N_\mu$.  Then as we saw earlier, $\xi^aC^c_{ab}\mu_c=0$, so
$\tilde{E}_a|_m(L_\xi)=0$.  Furthermore,
$\clift{\bar{E}_a}(\vlift{\bar{X}_i}(\R))=-C^b_{ac}p_bB^c_i$ as we showed
before, and $\Gamma=\xi^a\clift{\bar{E}_a}$ by assumption, so
\[
X_i|_m(L_\xi)=\Gamma(\vlift{\bar{X}_i}(\R))(m,0,\xi^a)=
\xi^a\clift{\bar{E}_a}(\vlift{\bar{X}_i}(\R))(m,0,\xi^a)
=-\xi^aC^b_{ac}\mu_bB^c_i(m,0,\xi^a)=0.
\]
Thus $m$ is a critical point of $L_\xi$.

Conversely, suppose that $m$ is a critical point of $L_\xi$.  Then
$C_{ab}^c\xi^b\mu_a=0$, so $\clift{\tilde{\xi}}$ is tangent to the
level set on which $(m,0,\xi^a)$ lies.  Furthermore, we have
$\Gamma(\vlift{\bar{X}_i}(\R))(m,0,\xi^a)=0$.  Recall that
$\Gamma=v^i\clift{\bar{X}_i}+\Gamma^i\vlift{\bar{X}_i}+v^a\clift{\bar{E}_a}$;
it follows that
\[
\Gamma(\vlift{\bar{X}_i}(\R))(m,0,\xi^a)=
\Gamma^j(m,0,\xi^a)\vlift{\bar{X}_i}(\vlift{\bar{X}_j}(\R))(m,0,\xi^a)+
\xi^a\clift{\bar{E}_a}(\vlift{\bar{X}_i}(\R))(m,0,\xi^a).
\]
But $\xi^a\clift{\bar{E}_a}(\vlift{\bar{X}_i}(\R))(m,0,\xi^a)
=-\xi^aC^b_{ac}\mu_bB^c_i(m,0,\xi^a)=0$, so
\[
\Gamma^j(m,0,\xi^a)\vlift{\bar{X}_i}(\vlift{\bar{X}_j}(\R))(m,0,\xi^a)=0.
\]
Since by assumption $\vlift{\bar{X}_i}\vlift{\bar{X}_j}(\R)$ is
non-singular, we have $\Gamma^i(m,0,\xi^a)=0$, and $(m,0,\xi^a)$ is a
relative equilibrium point.

If one is looking for relative equilibria with a given value of the
momentum $\mu$ it is appropriate to use the first method (searching
for critical points of the restriction of the energy function to the
level set $N_\mu$); if one is looking for relative equilibria with a
particular value of $\xi\in\g$ then the method described above is more
suitable.

\section{Some applications}

\subsection{Systems on Lie groups}

We now specialize to the case of an invariant Lagrangian system on a
Lie group $G$.  For such a system there are no conditions for relative
equilibria arising from the $X_i$, so the only condition for a point
$\tilde{\xi}_g\in TG$ to be a relative equilibrium point is that
$\xi^bC^c_{ab}\mu_c=0$ where $\mu$ is the value of the momentum at
$\tilde{\xi}_g$.  Thus the necessary and sufficient condition for a
relative equilibrium takes either of the following equivalent simple
forms:\ $\tilde{\xi}_g$ is a relative equilibrium point if and only if
the vector field $\clift{\tilde{\xi}}$ is tangent to the level set of
momentum in which the point $\tilde{\xi}_g$ lies, or equivalently if
and only if $\xi\in\g_\mu$, the algebra of the isotropy subgroup of
the momentum.

In the present case the fact that the Lagrangian is invariant means
that the dynamical system on $TG$ is determined by its reduction to
$T_eG\simeq\g$.  That is to say, the Euler-Lagrange equations can be
reduced to an equivalent set of equations on $\g$, the so-called
Euler-Poincar\'{e} equations \cite{MRbook}, which can be written
\[
\frac{d}{dt}\left(\fpd{l}{\xi^a}\right)=-C_{ab}^c\xi^b\fpd{l}{\xi^c}:
\]
here $l$ is the restriction of $L$ to $T_eG$, thought of as a function
on $\g$, and the $\xi^a$ here are the Cartesian coordinates on $\g$
determined by the basis $\{E_a\}$.  These equations, which are
first-order differential equations in the variables $\xi^a$, determine
in the regular case a vector field $\gamma$ on $\g$ from which the
Euler-Lagrange field $\Gamma$ on $TG$ can be reconstructed.  In fact a
curve $t\mapsto g(t)$ in $G$ is a base integral curve of $\Gamma$ if
and only if the curve $t\mapsto T\psi^M_{g(t)^{-1}}\dot{g}(t)$ in $\g$
is an integral curve of $\gamma$.

In this picture the relative equilibria are simply constant solutions
of the Euler-Poincar\'{e} equations, and these are points $\xi$ of
$\g$ at which
\[
C_{ab}^c\xi^b\fpd{l}{\xi^c}(\xi)=0.
\]
A solution of these equations determines a relative equilibrium
starting at $e$, or in other words a base integral curve of $\Gamma$
which coincides with a 1-parameter subgroup of $G$; but since
translates of relative equilibria are relative equilibria, this is
enough to give all relative equilibria.  Now
\[
p_c|_{T_eG}=\vlift{\tilde{E}_c}(L)|_{T_eG}=\fpd{l}{\xi^c},
\]
so the two approaches give the same results so far as relative
equilibria through the identity are concerned.

We discuss next the relations between the general criteria for finding
relative equilibrium points obtained earlier and the observations
above.  In order to do so we must first consider the identification of
$TG$ with $G\times\g$.  Since we are working with left actions the
fundamental vector fields are right, not left, invariant, so the use
of quasi-coordinates relative to a basis of fundamental vector fields
amounts to identifying $T_gG$ with $T_eG$ by right rather than left
translation.  On the other hand, when we say for example that $L$ is
invariant we mean that it is invariant under left translations.  Under
left translation, $\tilde{\xi}_g$ is identified with
$\ad_{g^{-1}}\xi$.  For any right-invariant function $F$ we have
$F(\tilde{\xi}_g)=F(\widetilde{\ad_{g^{-1}}\xi}|_e)$.  So if we denote
by $f$ the function on $\g$ obtained by restricting $F$ to $T_eG$ (and
identifying $T_eG$ with $\g$), then
$F(\tilde{\xi}_g)=f(\ad_{g^{-1}}\xi)$.

The energy $\E$ in this case is just
\[
\E(\tilde{\xi}_g)=\xi^a p_a(\tilde{\xi}_g)-L(\tilde{\xi}_g)
\]
(so $\E$ happens to coincide with $-\R$).  Now $\E$ is left-invariant,
and $\varepsilon$, its restriction to $\g$, is just
\[
\varepsilon(\xi)=\xi^a\fpd{l}{\xi^a}(\xi)-l(\xi).
\]
Notice that
\[
\fpd{\varepsilon}{\xi^a}(\xi)=\spd{l}{\xi^a}{\xi^b}(\xi)\xi^b
=\bar{g}_{ab}(\xi)\xi^b,
\]
where $\bar{g}_{ab}$ is the restriction of $g_{ab}$ to $T_eG\simeq\g$;
by assumption, the matrix $(\bar{g}_{ab})$ is non-singular everywhere
on $\g$.

The relative equilibrium points are the critical points of $\E^\mu$,
the restriction of $\E$ to the level set of momentum $N_\mu$.  To
express this result in terms of $\varepsilon$ we must determine those
points $(g,\xi)\in G\times\g\simeq TG$ which lie in $N_\mu$.  Now it
follows from the regularity assumptions that $N_\mu$ is (the image of)
a section of $TG\to G$, so that for each $g\in G$ there is a unique
$\xi\in\g$ such that $(g,\xi)\in N_\mu$.  It follows from equivariance
that $g$ and $\xi$ must satisfy $\ad_{g^{-1}}^* p(\tilde{\xi}_e)=\mu$,
or
\[
\fpd{l}{\xi^a}(\xi)=(\ad_{g}^*\mu)_a.
\]
This defines a map $G\to\g$, which is constant on left cosets of
$G_\mu$, the isotropy group of $\mu$ under the coadjoint action.  Let
$\g(\mu)\subset\g$ be the image of $G$ under this map.  Then the
relative equilibrium points in $T_eG$ with momentum $\mu$ are the
critical points of $\varepsilon$ restricted to $\g(\mu)$.

Now consider any curve in $N_\mu$, given in the form
$t\mapsto(g(t),\xi(t))$, such that $g(0)=e$; we set $\xi(0)=\xi_0$
and note that
\[
\fpd{l}{\xi^a}(\xi_0)=\mu_a.
\]
By differentiating the condition
\[
\fpd{l}{\xi^a}(\xi(t))=(\ad_{g(t)}^*\mu)_a
\]
with respect to $t$ and setting $t=0$ we obtain
\[
\bar{g}_{ab}(\xi_0)\dot{\xi}^b(0)=\eta^bC_{ba}^c\mu_c,
\]
where $\eta$ is the tangent vector to $t\mapsto g(t)$ at $t=0$,
considered as a point of $\g$.  We may choose $\eta$ arbitrarily, and
determine $\dot{\xi}(0)$ from this equation. The tangent vectors to
$\g(\mu)$ at $\xi_0$ are those of the form
\[
\bar{g}^{ac}(\xi_0)\eta^bC_{bc}^d\mu_d\vf{\xi^a}.
\]
It follows that $\xi_0$
will be a critical point of $\varepsilon|_{\g(\mu)}$ if and only if
\[
\bar{g}^{ac}(\xi_0)\eta^bC_{bc}^d\mu_d\fpd{\varepsilon}{\xi^a}(\xi_0)=0
\]
for all $\eta$. This gives back the same condition as before.

This approach is similar in spirit to that discussed by Arnold
\cite{Arn}, and indeed generalizes that approach insofar as the
finite-dimensional case is concerned since Arnold deals only with
kinetic energy Lagrangians defined by Riemannian metrics.

The locked Lagrangian for a system on a group $G$ is given by
$L_\xi(g)=L(\tilde{\xi}_g)$, for fixed $\xi$.  It follows from the
invariance assumption that $L_{\xi}(g)=l(\ad_{g^{-1}}\xi)$.  We can
think of the right-hand side as the restriction of $l$ to the orbit of
$\xi$ under the adjoint action of $G$ on $\g$, which we denote by
$G(\xi)$; that is, $L_\xi=l|_{G(\xi)}$.  According to Lewis's
criterion, $(g,\xi)$ is a relative equilibrium point if and only if
$g$ is a critical point of $L_\xi$.  Let $G_\xi$ be the isotropy group
of $\xi$ under the adjoint action; then $G(\xi)\simeq G/G_\xi$.
Clearly $L_\xi$ is constant on the fibres of the projection $\rho:G\to
G/G_\xi$, from which it follows that $g$ is a critical point of
$L_\xi$ if and only if $\rho(g)$ is a critical point of $l|_{G(\xi)}$.
Thus in this case Lewis's criterion can be restated in the following
form:\ $(g,\xi)$ is a relative equilibrium point if and only if
$\rho(g)$ is a critical point of $l|_{G(\xi)}$.  Lewis's criterion
again reduces to the condition $\xi^bC^c_{ab}\mu_c=0$, or more
succinctly $\langle[\eta,\xi],\mu\rangle=0$ for all $\eta\in\g$.  The
role of $G_\xi$ is revealed here by the observation that this
condition is automatically satisfied if $[\eta,\xi]=0$, that is, if
$\eta$ lies in the centralizer of $\xi$:\ but this is exactly the
algebra of $G_\xi$.

There is yet another way of arriving at the condition
$\xi^bC^c_{ab}\mu_c=0$.  The fundamental vector fields $\tilde{E}_a$
are the right translates of the $E_a$, considered as elements of
$T_eG$; they are not of course left-invariant.  We denote by
$\hat{E}_a$ the left translates of the $E_a$, which are
left-invariant.  The relation between these two sets of vector fields
on $G$ can be written $\hat{E}_a=A_a^b\tilde{E}_b$; the coefficients
are the matrix components of the adjoint map, and the condition of
invariance gives
\[
\tilde{E}_a(A_b^c)-C_{ab}^dA_d^c=0,
\]
where of course $A_b^c=\delta_b^c$ at $e$.
Now for any vector field $Y$ and function $f$ on a manifold $M$,
\[
\clift{(fX)}=f\clift{X}+\dot{f}\vlift{X},
\]
where $\dot{f}$ is the so-called total derivative of $f$, a function
on $TM$ given by
\[
\dot{f}=u^\alpha\fpd{f}{x^\alpha}=v^\alpha X_\alpha(f)
\]
for a vector field basis $\{X_\alpha\}$ with associated
quasi-velocities $v^\alpha$.
Thus
\[
\clift{\hat{E}_a}=A_a^b\clift{\tilde{E}_b}
+\xi^c\tilde{E}_c(A_a^b)\vlift{\tilde{E}_b}
=A_a^b\clift{\tilde{E}_b}+\xi^cC_{ca}^dA_d^b\vlift{\tilde{E}_b}.
\]
It follows that at the identity
\[
\clift{\hat{E}_a}(L)=\xi^cC_{ca}^bp_b.
\]
The necessary and sufficient conditions for $\xi$ to define a relative
equilibrium at $e$ may therefore be written
$\clift{\hat{E}_a}(L)(e,\xi)=0$.

We note in passing that if the Lagrangian is bi-invariant, that is,
invariant under both left and right translations, so that
$\clift{\hat{E}_a}(L)=0$ everywhere (as well as
$\clift{\tilde{E}_a}(L)=0$), then all base integral curves of
$\Gamma$ through $e$ coincide with 1-parameter subgroups, and
therefore all base integral curves are translates of 1-parameter
subgroups. These curves are just the geodesics of the canonical
torsionless connection on $G$, which is defined by
\[
\nabla_{\hat{E}_a}\hat{E}_b=\onehalf C^c_{ab}\hat{E}_c.
\]
The Euler-Poincar\'{e} equations reduce to
\[
C_{ab}^c\xi^b\fpd{l}{\xi^c}=0.
\]
Conversely, if $L$ is left-invariant ($\clift{\tilde{E}_a}(L)=0$) and
all base integral curves of its Euler-Lagrange field $\Gamma$ are
translates of 1-parameter subgroups then $L$ must be bi-invariant.
For it must certainly be the case that $\clift{\hat{E}_a}(L)(e,\xi)=0$
for all $\xi\in\g$.  But $\clift{\tilde{E}_b}\clift{\hat{E}_a}(L)=
\clift{\hat{E}_a}\clift{\tilde{E}_b}(L)=0$, so
$\clift{\hat{E}_a}(L)=0$ everywhere.  It is a well-known property of
invariant Riemannian metrics on a Lie group that the exponential map
determined by the Levi-Civita connection coincides with the
exponential in the group sense if and only if the metric is
bi-invariant.  The result above is a generalization of this property
to regular invariant Lagrangians.

The problem of the existence of relative equilibria for invariant
systems on Lie groups has been studied recently by several authors,
using differing terminology:\ Hern\'{a}ndez-Gardu\~{n}o {\em et al.}
\cite{HG} (for kinetic energy Lagrangians, i.e.\ geodesics of an
invariant Riemannian metric on a Lie group); Latifi \cite{Latifi} (for
invariant Finsler structures, under the name `homogeneous geodesics');
Szenthe \cite{Szenthe} (for a general invariant Lagrangian, under the
name `stationary geodesics').  Our results above incorporate the
particular cases in \cite{HG} and \cite{Latifi}.  Furthermore, our
results improve on those of Szenthe \cite{Szenthe} in that we do not
require one of the hypotheses, namely that the Lagrangian is a first
integral of its Euler-Lagrange field, in both his Proposition~2.2,
which (in different notation) gives the condition for a relative
equilibrium in the form $\clift{\hat{E}_a}(L)(e,\xi)=0$, and his
Theorem~2.3, which gives the condition in terms of critical points of
$l|_{G(\xi)}$.

\subsection{Simple mechanical systems}

A simple mechanical system is a Lagrangian system in which the
Lagrangian takes the familiar form $L=T-V$, where $T$ is the kinetic
energy associated with a Riemannian metric $g$ on $M$ and $V$ is the
potential energy, a function on $M$. Such a Lagrangian is necessarily
regular since its Hessian is effectively just the Riemannian metric.

In the case of a simple mechanical system we take as symmetry group
$G$ the group of diffeomorphisms of $M$ which are isometries of the
metric and leave the potential invariant.  We must assume of course
that $G$ acts freely and effectively on $M$.  We define the invariant
vector fields $X_i$ of a standard basis as follows.  The orthogonal
complements to the tangent spaces to the fibres of the principal
bundle $M\to B$ are the horizontal subspaces of a principal
connection, called the mechanical connection.  The $X_i$ are the
horizontal lifts to $M$, relative to the mechanical connection, of the
vector fields of some local basis on $B$.  We write
$g_{ab}=g(\tilde{E}_a,\tilde{E}_b)$, $g_{ij}=g(X_i,X_j)$; by
assumption, $g_{ai}=g(\tilde{E}_a,X_i)=0$.  Thus
\[
L(m,v)=\onehalf\left(g_{ij}(m)v^iv^j+g_{ab}(m)v^av^b\right)-V(m),
\]
where the $v$\,s are the quasi-velocities associated with the
standard basis, as before. It is clear that the $g_{ab}$ etc., which
are here defined as components of the metric, are also the
appropriate components of the Hessian of $L$. Since we
assume that $L$ is regular the matrix $(g_{ab}(m))$ is necessarily
non-singular in this case.

In the case of a simple mechanical system the components of momentum
are given simply by $p_a=g_{ab}v^b$.  The restriction of the Routhian
to a level set of momentum is
\[
\R^\mu=\onehalf g_{ij}v^iv^j-\left(V+\onehalf g^{ab}\mu_a\mu_b\right)
=\onehalf g_{ij}v^iv^j-V^\mu;
\]
$V^\mu$ is the so-called amended potential \cite{Simo}.
The restriction of the energy to a level set is given by
\[
\E^\mu=\onehalf g_{ij}v^iv^j+\onehalf g^{ab}\mu_a\mu_b+V =\onehalf
g_{ij}v^iv^j+V^\mu.
\]
The relative equilibrium points on the level set are determined by
the critical points of $\E^\mu$, and these are points of the form
$(m,0,\xi^a)$ where $\xi^a=g^{ab}\mu_b$ and $m$ is a critical point
of $V^\mu$. Now one of the conditions for a relative equilibrium
point is that $C^c_{ab}\xi^b\mu_c=C^c_{ad}g^{bd}(m)\mu_b\mu_c=0$;
this is in fact included in the condition for $m$ to be a critical
point of $V^\mu$. To see this, note that
\[
\tilde{E}_a(g^{bc})=-g^{bd}g^{ce}\tilde{E}_a(g_{de})
=g^{bd}g^{ce}\left(C_{ad}^fg_{ef}+C_{ae}^fg_{df}\right)
=g^{bd}C^c_{ad}+g^{ce}C_{ae}^b.
\]
It follows that
\[
\tilde{E}_a(\onehalf g^{bc}\mu_b\mu_c) =g^{bd}C^c_{ad}\mu_b\mu_c
\]
as required. So if $(m,0,\xi^a)$ is a relative equilibrium point on
the level set $p_a=\mu_a$, then $\mu_a=g_{ab}(m)\xi^b$, and $m$ must
be a critical point of the amended potential $V^\mu$.  Conversely,
if $m$ is a critical point of $V^\mu$ then $(m,0,\xi^a)$ is a
relative equilibrium point, where $\xi^a=g^{ab}(m)\mu_b$.

On the other hand, the locked Lagrangian $L_\xi$ is given by
\[
L_\xi=\onehalf g_{ab}\xi^a\xi^b-V;
\]
the quantity $V-\onehalf g_{ab}\xi^a\xi^b$ is the augmented or
effective potential \cite{Simo}, $V_\xi$. Then $(m,0,\xi^a)$ is a relative
equilibrium point if and only if $m$ is a critical point of $V_\xi$.
Notice that for any $w\in T_mM$, $w(g^{ab})=-g^{ac}g^{bd}w(g_{cd})$,
so that if $\xi^a=g^{ab}\mu_b$
\[
w(V^\mu)=w(V)-g^{ac}g^{bd}w(g_{cd})\mu_a\mu_b
=w(V)-w(g_{cd})\xi^c\xi^d=w(V_\xi),
\]
so the two criteria for the existence of a relative equilibrium
point are consistent.

Since $\mu_a=g_{ab}(m)\xi^b$, the condition $C^c_{ab}\xi^b\mu_c=0$ can
be written in the form $C^c_{ad}g^{bd}(m)\mu_b\mu_c=0$, as we have
already observed, and also in the form
$C^c_{ab}\xi^bg_{cd}(m)\xi^d=0$.  Now $g_{ab}$ may be regarded as
defining a function on $M$ taking its values in the space of symmetric
bilinear forms on $\g$, in the sense that for any $m\in M$,
$(g_{ab}(m))$ is the matrix of such a bilinear form with respect to
the basis $\{E_a\}$ of $\g$.  With this interpretation we can express
the condition $C^c_{ab}\xi^bg_{cd}(m)\xi^d=0$ equivalently as
$g(m)(\xi,[\xi,\eta])=0$ for all $\eta\in\g$.  This generalizes a
result of Szenthe's \cite{Szenthe1} for the case of an invariant
Riemannian metric on a Lie group, when this condition with $m=e$ is
the only condition for $\xi$ to determine a relative equilibrium
through the identity.

\subsection{Saari's conjecture}

We continue to discuss the case of a simple mechanical system.

The matrix-valued function $(g_{ab})$ is called the {\em locked
inertia tensor}.

It has been conjectured (see \cite{HG,Lawson}), on the basis of
certain results for the $N$-body problem, that `a Lagrangian simple
mechanical system with symmetry is at a point of relative equilibrium
if and only if the locked inertia tensor is constant along the
integral curve that passes through that point'. The original version
of this conjecture, in the context of the $N$-body problem, was
formulated by Saari; the version above is called the naive
generalization of Saari's conjecture.

It is evident from the formula
\[
\tilde{\xi}(g_{bc})=\xi^a\tilde{E}_a(g_{bc})
=-\left(\xi^aC_{ab}^dg_{cd}+\xi^aC_{ac}^dg_{bd}\right),
\]
which is part of Killing's equation for $\tilde{\xi}$, that if at a
relative equilibrium point $(m,0,\xi^a)$ we have $\xi^aC_{ab}^c=0$
(and not just $\xi^aC_{ab}^c\mu_c=0=C_{ab}^dg_{cd}\xi^a\xi^c$) then
the locked inertia tensor is constant along the corresponding integral
curve.  That is to say, if $(m,0,\xi^a)$ is a relative equilibrium
point for which $\xi$ belongs to the centre of $\g$ then the locked
inertia tensor is constant along the integral curve.  On the other
hand, our analysis suggests that it is unlikely that in general the
locked inertia tensor is necessarily constant along the integral curve
of a relative equilibrium.  So it seems unlikely that Saari's
conjecture holds in all generality; and indeed it is known to be
false.  In a refined version of Saari's conjecture formulated in
\cite{HG,Lawson} it is required only that $g_{ab}\xi^b$ is constant
along the integral curve.  It is clear that when $\xi$ does define a
relative equilibrium point, $g_{ab}\xi^b$ is constant along the
integral curve, because $g_{ab}\xi^b=\mu_a$ is the value of the
momentum at the relative equilibrium point and the integral curve
through the point lies in the same level set.  In fact, even in the
case of a general Lagrangian we have
\[
\clift{\tilde{\xi}}(g_{bc}\xi^c)
=-(\xi^aC_{ab}^dg_{cd}+\xi^aC_{ac}^dg_{bd})\xi^c
=-C_{ab}^dg_{cd}\xi^a\xi^c
\]
(using a formula which generalises the one at the beginning of this
paragraph), from which it is clear that $g_{bc}\xi^c$ is constant
along an integral curve of $\clift{\tilde{\xi}}$ if and only the
condition $C_{ab}^dg_{cd}\xi^a\xi^c=0$ holds.  This is indeed a
requirement for a point to be a relative equilibrium point, in the
case of a simple mechanical system; however, in general there is a
further requirement involving critical points of the augmented
potential.  But in the case of an invariant simple Lagrangian on a Lie
group $C_{ab}^dg_{cd}\xi^a\xi^c=0$ is the only condition for $\xi$ to
be a relative equilibrium point, so invariant simple Lagrangians on
Lie groups belong to the class of Lagrangian systems with symmetry for
which the refined Saari conjecture holds, as is pointed out in
\cite{HG}.

\subsubsection*{Acknowledgements}
The first author is a Guest Professor at Ghent University:\ he is
grateful to the Department of Mathematical Physics and Astronomy at
Ghent for its hospitality.

The second author is currently a Research Fellow at The University
of Michigan through a Marie Curie Fellowship. He is grateful to the
Department of Mathematics for its hospitality. He also acknowledges
a research grant (Krediet aan Navorsers) from the Fund for
Scientific Research - Flanders (FWO-Vlaanderen), where he is an
Honorary Postdoctoral Fellow.

\end{document}